\documentclass{article}
\usepackage[utf8]{inputenc}
\usepackage{amssymb,amsmath,amsthm}
\theoremstyle{definition}

\theoremstyle{theorem}
\newtheorem{theorem}{Theorem}

\newtheorem{conjecture}[theorem]{Conjecture}

\DeclareMathOperator{\Prob}{Prob}

\title{Note: a counterexample to a conjecture of Gilmer which would imply the union-closed conjecture}
\author{David Ellis\footnote{School of Mathematics, University of Bristol, UK. Email: \texttt{david.ellis@bristol.ac.uk}}}
\date{22nd November 2022}

\begin{document}

\maketitle

\begin{abstract}
    In this very short note, we give a counterexample to a recent conjecture of Gilmer which would have implied the union-closed conjecture.
\end{abstract}
\section{Introduction}
We say a family $\mathcal{F}$ of sets is {\em union-closed} if $A \cup B \in \mathcal{F}$ whenever $A,B \in \mathcal{F}$. The celebrated Union-Closed conjecture of Frankl states that if $n \in \mathbb{N}$ and $\mathcal{F} \neq \{\emptyset\}$ is a union-closed family of subsets of $\{1,2,\ldots,n\}$, then there exists $i \in \{1,2,\ldots,n\}$ such that at least half of the sets in $\mathcal{F}$ contain $i$. Gilmer \cite{gilmer} recently obtained a breakthrough on this conjecture, proving that for any union-closed family $\mathcal{F}$ of subsets of $\{1,2,\ldots,n\}$ with $\mathcal{F} \neq \{\emptyset\}$, there exists an element $i$ contained in at least $0.01|\mathcal{F}|$ of the sets in $\mathcal{F}$. (We refer
the reader to the survey of Bruhn and Schaudt \cite{bs} for
survey of work prior to Gilmer's, on the problem.)

Gilmer's proof is an elegant entropy argument. At the end of his paper, Gilmer makes the following (information-theoretic) conjecture which would immediately imply the union-closed conjecture. 

\begin{conjecture}[Gilmer]
Let $A, B$ be i.i.d.\ samples from a distribution over a family of
subsets of $\{1,2,\ldots,n\}$. Assume that $\Prob[i \in A] < \tfrac{1}{2}$ for all $i$, and that $H(A) > 0$. Then
$$H(A \cup B) + D(A \cup B||A) > H(A).$$
\end{conjecture}
Here, $H(A)$ denotes the entropy of $A$; recall that for a probability distribution $p = (p_x)_{x \in X}$ over a finite set $X$, the entropy of $p$ is defined by
$$H(p) = \sum_{x \in X} p_x \log_2(1/p_x).$$
Also, $D(q||p)$ denotes the {\em Kullback-Leibler divergence} of $q$ from $p$, for probability distributions $(p_x)_{x \in X}$ and $(q_x)_{x \in X}$ over a finite set $X$; recall that this is defined by
$$D(q||p) = \sum_{x \in X}q_x \log_2 (q_x/p_x).$$
The purpose of this very short note is to give a (simple) counterexample to this conjecture (with $n=2$). Another counterexample (for large $n$) was independently and simultaneously obtained by Sawin.

\section{The counterexample}
As usual, we write $[n]: = \{1,2,\ldots,n\}$, and for a set $S$, we write $\mathcal{P}(S)$ for the power-set of $S$. We first note that, writing $p$ for the distribution of $A$ in Gilmer's conjecture, and $q$ for the distribution of $A \cup B$, the left-hand side of the conjectured inequality is equal to
$$\sum_{x \subset [n]}q_x \log_2(1/q_x)+\sum_{x \subset [n]}q_x \log_2(q_x/p_x) = \sum_{x \subset [n]} q_x \log_2(1/p_x),$$
and therefore Gilmer's conjecture is equivalent to
\begin{equation} \label{eq:quantity} \sum_{x \subset [n]} q_x \log_2(1/p_x) - \sum_{x \subset [n]}p_x \log_2(1/p_x) >0.\end{equation}
We first give a probability distribution $p$ on $\mathcal{P}([2])$ such that, if $A$ and $B$ are i.i.d.\ samples from $p$, then $\Prob[1 \in A] = \Prob[2 \in A] = \tfrac{1}{2}$, but writing $q$ for the distribution of $A \cup B$, the quantity on the left-hand side of (\ref{eq:quantity}) satisfies
\begin{equation}\label{eq:quantity2} \sum_{x \subset [n]} q_x \log_2(1/p_x) - \sum_{x \subset [n]}p_x \log_2(1/p_x) < -0.04.\end{equation}
The distribution $p$ is as follows:
\begin{equation}\label{eq:p-defn} p(\emptyset) = p(\{1,2\}) = x,\quad p(\{1\}) = p(\{2\}) = \tfrac{1}{2}-x,\end{equation}
where (for concreteness) we take $x=0.3$. (We use the variable $x$ to make the construction more readable.)

It remains to observe that an arbitrarily small perturbation ($p'$, say) of the probability distribution $p$ satisfies the hypotheses of Gilmer's conjecture and yet also has the above quantity (\ref{eq:quantity2}) being negative; indeed, we may replace $p$ by the distribution $p'$ defined by
$$p'(\emptyset) = x,\quad p'(\{1,2\}) = x-2\epsilon,\quad p'(\{1\}) = p'(\{2\}) = \tfrac{1}{2}+\epsilon-x,$$
for $\epsilon$ a sufficiently small positive number.

We now check that we do indeed have
$$\sum_{x \subset [2]} q_x \log_2(1/p_x) - \sum_{x \subset [2]}p_x \log_2(1/p_x)<-0.04,$$
for the probability distribution $p$ defined in (\ref{eq:p-defn}) above. Indeed, we have
$$q_{\emptyset} = \Prob[A = B = \emptyset] = x^2,$$
and
\begin{align*} q_{\{1\}}& = \Prob[A =\emptyset \text{ and }B = \{1\}]+\Prob[A = \{1\} \text{ and }B = \emptyset]+\Prob[A=B=\{1\}]\\
&= 2x(\tfrac{1}{2}-x)+(\tfrac{1}{2}-x)^2\\
&= x-2x^2 + \tfrac{1}{4}-x+x^2\\
& = \tfrac{1}{4}-x^2,
\end{align*}
and by symmetry $q_{\{2\}} = q_{\{1\}} = \tfrac{1}{4}-x^2$. Since $q$ is a probability distribution, we have
$$q_{\{1,2\}} = 1-q_{\emptyset}-q_{\{1\}}-q_{\{2\}} = 1-x^2-2(\tfrac{1}{4}-x^2) = \tfrac{1}{2}+x^2.$$

Substituting in the above values, we have
\begin{align*} &\sum_{x \subset [2]} q_x \log_2(1/p_x) - \sum_{x \subset [2]}p_x \log_2(1/p_x)\\
& = (\tfrac{1}{2}+2x^2-2x)\log_2(1/x) + (-\tfrac{1}{2}-2x^2+2x) \log_2(1/(\tfrac{1}{2}-x)),\end{align*}
which is indeed less than -0.04 when $x=0.3$, as claimed.

\end{document}